\documentclass[10pt,a4paper]{amsart}
\usepackage[T1]{fontenc}
\usepackage{lmodern}
\usepackage[margin=27mm,top=24mm,bottom=26mm]{geometry}
\usepackage{amsmath,amssymb,amsthm,mathtools}
\usepackage{xcolor}
\definecolor{citation}{rgb}{0.2,0.55,0.2}
\definecolor{formula}{rgb}{0.1,0.2,0.6}
\definecolor{urlc}{rgb}{0,0,0.45}
\usepackage[colorlinks=true,linkcolor=formula,citecolor=citation,urlcolor=urlc]{hyperref}

\newtheorem{theorem}{Theorem}[section]
\newtheorem{proposition}[theorem]{Proposition}
\newtheorem{lemma}[theorem]{Lemma}
\numberwithin{equation}{section}
\AddToHook{env/abstract/begin}{\hyphenpenalty=10000\relax}
\newcommand{\R}{\mathbb R}\newcommand{\E}{\mathbb E}\newcommand{\Sph}{\mathbb S}
\newcommand{\one}{\mathbf 1}\newcommand{\dd}{\,\mathrm d}
\DeclareMathOperator{\sgn}{sgn}
\allowdisplaybreaks[2]\emergencystretch=2em

\title[A singular fractional perimeter minimizer in dimension seven]{A singular Lawson cone minimizing the fractional perimeter in dimension seven}

\author{Giampiero Palatucci}
\address{Department of Mathematical, Physical and Computer Sciences, University of Parma, Parco Area delle Scienze 53/A, 43124 Parma, Italy}
\email{giampiero.palatucci@unipr.it}

\subjclass[2020]{Primary 49Q05. Secondary 49Q20, 53A10}
\keywords{Fractional perimeter, Lawson cones, singular minimizers, nonlocal calibrations}
\thanks{The author acknowledges the support of the INdAM GNAMPA group.}
\date{23 November 2026}

\begin{document}

\begin{abstract}
We prove that the nonlocal Lawson cone constructed by D\'avila, del Pino,
and Wei ({\it J.~Differential~Geom.}~2018) in $N=7$ minimizes the fractional perimeter for all sufficiently small orders~$s$.
Thus nonlocal minimal surfaces can have singularities already in dimension~$7$. This contrasts both with the classical theory, where the first singular example occurs in dimension~$8$ (Bombieri--De~Giorgi--Giusti, {\it Invent.~Math.}~1969), and with the fractional theory for orders~$s$ close to~$1$ (Caffarelli--Valdinoci, {\it Adv.~Math.}~2013).
\end{abstract}
\maketitle
\medskip\medskip\section{Introduction}
We prove that a nonlocal Lawson cone in $\R^7$, constructed by D\'avila, del Pino and Wei in~\cite{DDW}, minimizes the fractional perimeter for all sufficiently small orders~$s$.

\medskip

In their foundational paper~\cite{CRS}, Caffarelli, Roquejoffre and Savin developed a regularity
theory for sets that minimize the fractional perimeter. The relative perimeter $P_s(E;\Omega)$
is obtained by integrating $|x-y|^{-N-s}$ over pairs $(x,y)\in E\times E^c$ with at least one
point in $\Omega$. In the corresponding Plateau problem the datum is prescribed on the whole
complement of $\Omega$, since even points far outside $\Omega$ enter the energy. What is the
least ambient dimension in which a minimizing boundary can have an interior singularity? To
see why cones enter the question, consider a minimizer on successively smaller scales about an
interior boundary point. Every blow-up limit is a minimizing cone, and if one of them is a
half-space, the original boundary is regular near the point; see in particular Theorems~9.2 and~9.4 there.
Conversely, a nonflat minimizing cone already supplies an example, since its vertex is
singular. We refer to \cite[Chapters~V--VI]{FAG} for an account of the relative perimeter and
its exterior data.

\medskip

When $s$ is close to one, the main contribution to the rescaled perimeter comes from pairs of
nearby points. This is familiar from the Gagliardo seminorm; see~\cite{DNPV12}. For $u\in C_c^\infty(\R^N)$ and
$1<p<\infty$, the limit of $(1-s)[u]_{W^{s,p}(\R^N)}^p$ is the integral of $|\nabla u|^p$ up
to a constant depending on $N$ and $p$, as seen in~\cite{BBM}. For sets, the precise variational statement
on a bounded Lipschitz domain is that $(1-s)P_s(\,\cdot\,;\Omega)$ $\Gamma$-converges to the
classical perimeter in $\Omega$ times the volume of the unit ball in $\R^{N-1}$; see
\cite[Theorem~2]{ADPM}. If the exterior set is held fixed and has locally finite perimeter,
the limit also counts the portion of $\partial\Omega$ on which the inner and outer traces
disagree; see~\cite{DLPPP}.

\smallskip

The geometric step is to turn convergence of the energies into regularity of the minimizing
boundaries. This requires control over the scales at which a boundary becomes flat. Caffarelli
and Valdinoci obtained an improvement of flatness uniform as $s\to1^-$. Together with
compactness, it implies that in ambient dimension $N\le7$ every nontrivial minimizing cone is
a half-space and every minimizing boundary is regular in the interior for $s$ sufficiently
close to one; this is done in~\cite[Theorems~2 and~3]{CV}. This recovers the classical regularity range. For
the classical perimeter the bound is sharp. Simons' Theorem excludes a singular minimizing
cone through dimension seven; see \cite[Sections~6.1--6.2]{Simons}, and Bombieri, De Giorgi and
Giusti in~\cite{BDGG} proved that the Simons cone in $\R^4\times\R^4$ is area minimizing, giving
the first singular example in dimension eight.

\smallskip

Novaga and Paolini~\cite[Theorem~1.1]{NP} recently proved that either phase of the same
cone minimizes the fractional perimeter under compact perturbations for all $s$ sufficiently
close to one. Together with the regularity theorem of Caffarelli and Valdinoci, this shows
that for such $s$ the least ambient dimension in which a minimizing boundary can have an
interior singularity is eight.

The flatness conclusion in dimensions $3\le N\le7$ also holds under the weaker assumption
of stability. For stationary cones with $C^2$ boundary away from the vertex, Cabr\'e, Cinti
and Serra proved this in~$\R^3$; see \cite[Theorem~1.2]{CCS}. Chan, Dipierro, Serra and Valdinoci extended the
result to $\R^4$ in~\cite{CDSV}, and the recent work of Dipierro, Serra, Tsiamis and Valdinoci
covers the entire range $3\le N\le7$ in~\cite[Theorem~1.1]{DSTV}.

\medskip

Much less is known for the full range $0<s<1$. In the plane, Savin and Valdinoci proved that
every nontrivial minimizing cone is a half-plane, whatever the order. By dimension reduction
this gives an upper bound of $N-3$ for the Hausdorff dimension of the interior singular set in
every $N\ge3$; see \cite[Theorem~1 and Corollary~2]{SV}. The planar classification has a
consequence in three dimensions as well. Figalli and Valdinoci used it to show that every
entire minimizing graph in $\R^3$ is a plane; see \cite[Corollary~1.3]{FV}. These results leave
open what can happen to a general minimizing boundary in higher dimensions as the order
decreases.

\medskip

As $s\to0^+$, the interaction with distant points becomes decisive. Indeed, the radial
mass of the kernel beyond any fixed radius $R>1$ is $\int_R^\infty r^{-1-s}\dd r=R^{-s}/s$.
Even this part survives multiplication by $s$. In the whole space, Maz'ya and Shaposhnikova
in \cite{MS} showed that for $u\in C_c^\infty(\R^N)$ and $1<p<\infty$, $s[u]_{W^{s,p}(\R^N)}^p$ tends to a
multiple of $\|u\|_{L^p(\R^N)}^p$. For a relative perimeter, the exterior datum also
enters this limit. What matters at leading order is how much of the space at infinity each
phase occupies. To state the formula in dimension seven, write $\omega_6=\mathcal H^6(\Sph^6)$
and $B_1=B_1(0)$. When the limit exists, the normalized mass of $E$ at infinity is
\[
 \theta
       =\lim_{s\to0^+}
     \frac{s}{\omega_6}\operatorname{Tail}_{s/2,2}(\one_E;0,1)=\lim_{s\to0^+}\frac{s}{\omega_6}
       \int_{E\setminus B_1}\frac{\dd y}{|y|^{7+s}}\,,
\]
where $\operatorname{Tail}_{\sigma,p}$ is the nonlocal tail introduced in \cite{DKP14,DKP16}.
With $p=2$ and $\sigma=s/2$, the tail of $u$ at the unit ball is the integral of
$|u(y)|\,|y|^{-7-s}$ over $|y|>1$, so the tail of the constant one is $\omega_6/s$. Changing~$E$ in a bounded region leaves $\theta$ unchanged. For a bounded $C^{1,\gamma}$ domain
$\Omega\subset\R^7$, with $0<\gamma<1$, and~$P_{s_0}(E;\Omega)<\infty$ for some $s_0\in(0,1)$,
the limit computed in \cite[Theorem~2.5]{DFPV} is
\[
 \lim_{s\to0^+}\frac{sP_s(E;\Omega)}{\omega_6}
   =(1-\theta)|E\cap\Omega|+\theta|\Omega\setminus E|.
\]

This formula singles out $\theta=1/2$ when one seeks a minimizing interface that survives in
the limit. Indeed, with exterior data fixed and no volume constraint, the right-hand side is
minimized by the empty phase in $\Omega$ if $\theta<1/2$, and by the full phase if
$\theta>1/2$. Only when $\theta=1/2$ can both phases be present in a minimizer of the limiting
functional. In that case the value is $|\Omega|/2$ for every interior configuration, and so
the leading term says nothing about the shape of the interface. The first-order
$\Gamma$-expansion in \cite{DLPPP} describes where this information enters. After the
finite-$s$ contribution of the tail to the volume term is subtracted, the limit retains both
the interaction across the interface inside $\Omega$ and its interaction with the prescribed
exterior. If the normalized tail has a finite first-order correction about $\theta$, this
correction contributes an explicit volume term in the expansion about the limiting mass.

\smallskip

For cones, the balance has a simple geometric meaning. Integration in polar coordinates shows
that the normalized mass of a cone at infinity is the fraction of the unit sphere occupied by
its phase. Tsiamis in~\cite[Theorem~3.5]{Tsiamis} proved a corresponding compactness theorem in $\R^N$, $N\ge3$. Given a
sequence of stationary stable cones with nonempty, closed, embedded $C^2$ links and orders
tending to zero, the indicators of their spherical phases converge along a subsequence in
$L^1$ to that of a set occupying half the sphere.

Stability is particularly restrictive in low dimensions in this regime. Caselli proved that in
the plane, for small $s$, a nontrivial stationary cone stable away from the vertex is a
half-plane; \cite{Caselli}. Tsiamis recently proved that in $\R^3$, for all sufficiently small
$s$, every nontrivial locally minimizing cone is a half-space. The same conclusion holds for
stationary stable cones with nonempty $C^2$ boundary off the vertex;
\cite[Theorem~1.1]{Tsiamis}. The dimension enters when the stability inequality is passed to
the spherical limit. On~$\Sph^2$ it forces the limiting link to be a great circle, and a
rigidity estimate near the equator then gives flatness for positive orders.

\medskip

In higher dimensions it remained open whether a fractional perimeter minimizer could develop a
singularity below the classical threshold. Dimension seven had already emerged, however, in
the work of D\'avila, del Pino and Wei \cite[Theorems~3 and~4]{DDW}. For each splitting
$\R^m\times\R^n$, $m,n\ge1$, and each $0<s<1$, they found a unique aperture for which the
corresponding Lawson cone is stationary. Their second-variation analysis revealed a change
precisely at $m+n=7$. For sufficiently small orders these cones are unstable in dimensions at
most six and stable in dimension seven. This gave a candidate for a singular minimizer. The
second variation tests smooth perturbations, and it remained to decide whether an arbitrary
change inside a ball could lower the perimeter.

We prove here that the cone belonging to the splitting $\R^3\times\R^4$ does minimize the
fractional perimeter. Write $x=(Y,Z)$ with $Y\in\R^3$ and $Z\in\R^4$, corresponding to
$(m,n)=(4,3)$ in \cite{DDW}. If $\alpha_s$ is the stationary slope, the cone has equation
$|Y|=\alpha_s|Z|$. It will be convenient to use its latitude
$b_s=\alpha_s/\sqrt{1+\alpha_s^2}$ on the unit sphere. With $B_s=b_s^2$, write
\[
 Q_{B_s}(Y,Z)=B_s|Z|^2-(1-B_s)|Y|^2,\qquad
 C_s=\{Q_{B_s}>0\}=\{(Y,Z):|Y|<\alpha_s|Z|\}.
\]
For a ball $D$ and a measurable set $E$ we use the relative perimeter
\begin{equation}\label{eq:defper}
 P_s(E;D)=\frac12
 \iint_{\R^{14}\setminus(D^c\times D^c)}
 \frac{|\one_E(x)-\one_E(y)|}{|x-y|^{7+s}}\,dx\,dy .
\end{equation}

\begin{theorem}\label{th:main}
There exists $s_*\in(0,1)$ such that for every $0<s<s_*$, every ball $D\subset\R^7$ and every measurable $E$ agreeing with $C_s$ almost everywhere in $D^c$,
\[
 P_s(C_s;D)\le P_s(E;D).
\]
The same assertion holds for the complementary phase.
\end{theorem}

The singularity in this example is confined to the vertex. Indeed, $\nabla Q_{B_s}$ is nonzero at
every other point of the boundary. At the origin every dilation gives back the cone itself,
and its link $b_s\Sph^2\times\sqrt{1-b_s^2}\,\Sph^3$ spans $\R^7$, so there can be no tangent
hyperplane. 

\subsection*{Idea of the proof}
The comparison in Theorem~\ref{th:main} must hold even for competitors that change the
topology of the interface. We obtain it by placing the cone in a family of level surfaces and
computing the sign of their fractional curvature. This is the level-set calibration of Cabr\'e
\cite[Section~2]{Cabre}. The classical quartic $|U|^4-|V|^4$ used for the Simons cone
\cite{DPP,MM,Mir} also appears in the proof of Novaga and Paolini~\cite[Theorem~1.2]{NP}.
For $s$ close to one they use its positive levels as strict one-sided barriers and obtain
the barriers on the other side by interchanging the two $\R^4$ factors. The factorization
$|U|^4-|V|^4=(|U|^2+|V|^2)(|U|^2-|V|^2)$, together with the foliations and quantitative
calibrations for classical Lawson cones \cite{Dav,DPM}, suggests taking $G_{B_s}(x)=|x|^2Q_{B_s}(x)$
for the $3+4$ cone. The zero level of $G_{B_s}$ is the same as that of $Q_{B_s}$, while along each
ray its nonzero values grow as the fourth power of the radius. We shall show that its leaves
have negative signed curvature in $C_s$
and positive signed curvature in $\{Q_{B_s}<0\}$.

On the unit sphere the curvature depends only on the latitude $v=|Y|^2$, by invariance under
rotations of each block. As $s\to0^+$, the stationarity equation forces $b_s^2$ to
approach the latitude $b^2$ that divides $\Sph^6$ into two regions of equal area. This is the
balance $\theta=1/2$ for our limiting cone. In Section~\ref{sec:reduction} we subtract the conical tail before taking the limit. The resulting curvature and its derivative in the latitude near the link are continuous in both the aperture and the order, up to $s=0$. We compare the renormalized quartic curvature with its value on the link. At the stationary aperture this difference recovers the curvature for $s>0$. If the limiting difference has the required sign on each side and a positive derivative in $v$ at the link, continuity carries the sign to every sufficiently small positive order.

For the quadratic curvature we pair opposite rays through the observation point. Subtracting the value at the link cancels the constant terms before division by $s$ and leaves a logarithmic integral as $s$ tends to zero (Lemma~\ref{lem:four-curvature}).

To compare the curvatures, choose the common sign of the quadratic and quartic functions so
that they are positive in the phase under consideration. The gain is the quadratic curvature
minus the quartic curvature. It can be seen by looking at the spheres of radii $r$ and $1/r$
together. At an observation point on the unit sphere the two functions have the same value. In
the curvature difference, the contributions from inside and outside the sphere have opposite
signs. After inversion of the radius, a reflection pairs their intervals of integration in
latitude. With the radial weights included, the outgoing contribution dominates its reflected
partner. When the reflected latitude falls beyond an axis, the outgoing contribution is
retained in full (Proposition~\ref{prop:quarticmirror}).

It remains to show that the gain exceeds the quadratic curvature. At the limiting aperture we divide by the distance to the link before approaching it. An explicit radial integral then shows that the excess, so normalized, has a strictly positive lower limit. Farther out on the three-dimensional block side, the quadratic curvature divided by that distance has a decreasing upper bound. On the four-dimensional block side the comparison is more delicate. For each angular direction, the radius at which the reflected point disappears is also the endpoint of the integral that measures the loss. The outgoing gain and the reflection loss near the link can thus be estimated in terms of the same moment of this radius. Towards the axis a concavity argument reduces the last inequality to its two endpoint values (Sections~\ref{sec:prelim} and~\ref{sec:meridians}). The estimates cover both sides up to the axes. The passage to positive orders and the level-set calibration then give the perimeter comparison in Theorem~\ref{th:main}.

\medskip\section{Reduction to the limiting aperture}\label{sec:reduction}
Let $v(x)=|Y|^2/|x|^2$, $B_s=\alpha_s^2/(1+\alpha_s^2)$ and
\[
 F_{d,A}(x)=|x|^d(A-v(x)),\qquad d\in\{2,4\},\quad 0<A<1.
\]
Set $F_{d,A}(0)=0$, $Q_A=F_{2,A}$ and $G_A=F_{4,A}$. Thus $C_s=\{G_{B_s}>0\}$ up to a null set. Let $\omega_6=\mathcal H^6(\Sph^6)$, $d\sigma=d\mathcal H^6/\omega_6$ on $\Sph^6$, and
\[
 H_s[F](x)=\omega_6^{-1}\operatorname{pv}\!\int_{\R^7}
 \frac{\sgn(F(x)-F(y))}{|x-y|^{7+s}}\,dy .
\]
For homogeneous $F$ of positive degree, $H_s[F](x)=|x|^{-s}H_s[F](x/|x|)$.

\subsection*{The calibration}
With this orientation it suffices to show
\begin{equation}\label{eq:needed}
 (v(\xi)-B_s)\,H_s[G_{B_s}](\xi)\ge0\qquad (\xi\in\Sph^6).
\end{equation}
Indeed, set $\Phi(x,y)=\sgn(G_{B_s}(x)-G_{B_s}(y))$ and $\delta=\one_E-\one_{C_s}$.
For any pair of points,
$|\one_E(x)-\one_E(y)|\ge\Phi(x,y)(\one_E(x)-\one_E(y))$,
with equality almost everywhere for $C_s$. Truncate the kernel in \eqref{eq:defper} to $|x-y|>\epsilon$ and subtract the equality from the inequality. Antisymmetry of $\Phi$ and the fact that $\delta=0$ on $D^c$ give, after removal of the truncation,
\begin{equation}\label{eq:cal}
 P_s(E;D)-P_s(C_s;D)\ge
 \omega_6\int_D\delta(x)H_s[G_{B_s}](x)\,dx .
\end{equation}
Here the cone has finite relative perimeter. For $|h|<1$ a change of phase between $x\in D$ and $x+h$ requires
$\big||Y|-\alpha_s|Z|\big|\le (1+\alpha_s)|h|$; integration first in $|Y|$ bounds the measure of this set by $O_D(|h|)$. The interaction with $|h|>1$ is integrable for $s>0$.
On the unit sphere $\nabla G_{B_s}$ never vanishes, since its radial component is $4(B_s-v)$ and $|\nabla_{\Sph}v|=2\sqrt{v(1-v)}$. On a sphere of radius $\rho$ about $\xi$, Taylor's formula confines disagreement with the tangent half-space to an angular band of measure $O(\rho)$, uniformly in $\xi$. The half-space has zero signed mean. Integrating $O(\rho^{-s})$ in $\rho$ and scaling yields the cutoff-independent majorant $K_s|x|^{-s}\in L^1(D)$ for the truncated curvatures. Dominated convergence on the right and monotone convergence for each perimeter give \eqref{eq:cal}, also in the extended sense for competitors of infinite perimeter. Under \eqref{eq:needed} the integrand on the right is nonnegative. Complementation gives the other phase.

\subsection*{The balanced latitude}
On a fixed meridian take $\xi(v)=(\sqrt v\,e_1,\sqrt{1-v}\,e'_1)$.
For a uniform $\eta\in\Sph^6$, the source latitude $a=v(\eta)$ and the first coordinates $U,V$ of its unit block directions are independent. Their densities are
\begin{equation}\label{eq:laws}
 w(a)=\tfrac{15}{4}\sqrt a(1-a)\quad(0<a<1),\qquad
 \tfrac12,\quad \tfrac2\pi\sqrt{1-t^2}\quad(-1<t<1),
\end{equation}
respectively. This follows by normalizing a standard Gaussian vector, whose two squared block lengths have independent gamma laws with parameters $3/2$ and $2$. Thus $\xi(v)\cdot\eta=\sqrt{va}\,U+\sqrt{(1-v)(1-a)}\,V$. We extend $w$ by zero and write $\E$ for expectation.

Let $\mathsf M_A=\operatorname{diag}(-(1-A)I_3,A I_4)$. On the rays $\xi\pm\rho\eta$, the two signed contributions to the curvature of $Q_A$ cancel until $\rho=2|\mathsf M_A\xi\cdot\eta|/|Q_A(\eta)|$, after which their mean is $-\sgn Q_A(\eta)$. Radial integration gives
\begin{equation}\label{eq:ray}
 H_s[Q_A](\xi)=-\frac1s\int_{\Sph^6}\sgn Q_A(\eta)
        \left|\frac{Q_A(\eta)}{2\mathsf M_A\xi\cdot\eta}\right|^s d\sigma(\eta).
\end{equation}
The inverse $s$-moment of a spherical coordinate is finite for $s<1$, justifying removal of the centered cutoff.

At $A=B_s$ and $v=B_s$, \eqref{eq:ray} vanishes by stationarity. On any sequence $s\to0^+$ choose a subsequence on which $B_s\to A_0\in[0,1]$ and the unit normals $n_s=\mathsf M_{B_s}\xi(B_s)/|\mathsf M_{B_s}\xi(B_s)|$ converge. After cancellation of $(2|\mathsf M_{B_s}\xi(B_s)|)^{-s}$ the signed integrand converges almost everywhere to $\sgn(A_0-a)$. It is bounded in $L^2(d\sigma)$ for $s\le1/4$, since $|Q_{B_s}(\eta)|\le1$ and the distribution of $n_s\cdot\eta$ is independent of $s$. Hence $\int_0^{A_0}w=1/2$. The unique median is interior, so $B_s\to B=b^2$, where
\begin{equation}\label{eq:median}
 5b^3-3b^5=1,\qquad C=1-B,\qquad d_*=\tfrac32 bC.
\end{equation}
The polynomial on the left increases on $(0,1)$. Evaluation at $9/14$ and $\sqrt{5/12}$ gives
\begin{equation}\label{eq:apbounds}
 9/14<b<\sqrt{5/12}<13/20,\qquad 7/10<B/C<5/7,\qquad d_*<17/30 .
\end{equation}
For the last bound use that $b(1-b^2)$ decreases for $b>1/\sqrt3$.

We shall also use the latitude density bound
\begin{equation}\label{eq:weight}
 w(a)\ge w(B)>45/32\qquad (2B/3\le a\le B).
\end{equation}
Indeed $w$ is concave and $w(2B/3)/w(B)=\sqrt{2/3}(1+B/(3C))>1$ by $B/C>7/10$. The lower bound at $B$ follows from $w(B)=15bC/4$, $b>9/14$ and $C>7/12$.

\subsection*{Removing the conical tail}
Write $S_d=\sgn(F_{d,A}(\xi(v))-F_{d,A}(y))$ and $S_\infty=\sgn(v(y)-A)$.
For $0\le s\le1/2$ set
\[
 R_d(s,A,v)=\omega_6^{-1}\operatorname{pv}\!\int_{\R^7}
 \left(\frac{S_d}{|y-\xi(v)|^{7+s}}
       -\one_{\{|y|>2\}}\frac{S_\infty}{|y|^{7+s}}\right)dy.
\]
For $s>0$ the term subtracted from the curvature is independent of $v$ and $d$.
At $v=A$ the two level sets have the same positive phase, so $R_2(s,A,A)=R_4(s,A,A)$.

\begin{lemma}\label{lem:c1}
For $d=2,4$, $R_d$ is continuous in $(s,A,v)$ on
$[0,1/2]\times(0,1)\times[0,1]$. For $0<v<1$ its derivative in $v$ exists and is jointly continuous up to $s=0$.
\end{lemma}
\begin{proof}
Choose a smooth radial cutoff $\psi(|y-\xi|)$ equal to one near $\xi$ and supported in a small ball about it.
In the local part write $y=\xi+\rho\eta$ and choose coordinates depending smoothly on the parameters in which the first axis is $\nabla F_{d,A}(\xi)/|\nabla F_{d,A}(\xi)|$. The radial and spherical components of the gradient are $d(A-v)$ and $-\nabla_{\Sph}v$, so this choice is possible locally at every $\xi$, including the axes.
In the coordinates $\eta=(t,\sqrt{1-t^2}\,\zeta)$ of this frame the normalized level difference is
\[
 \frac{F_{d,A}(\xi+\rho\eta)-F_{d,A}(\xi)}
      {\rho|\nabla F_{d,A}(\xi)|}=t+\rho L(\rho,t,\zeta,A,\xi).
\]
Taylor's Formula gives a bounded $L$, with bounded first derivatives for $|t|\le1/2$ on compact parameter sets. For small $\rho$ this expression has a unique zero $t_0$ in $t$ and $|t_0|=O(\rho)$. Implicit differentiation gives $\partial_v t_0=O(\rho)$ for $v$ in compact subintervals of $(0,1)$. After subtraction of the tangent half-space, the angular mean is, up to sign, twice the mean in $\zeta$ of $\int_0^{t_0} c_7(1-t^2)^2\,dt$, where $c_7(1-t^2)^2$, with $c_7=15/16$, is the density of a coordinate on $\Sph^6$. Thus both the local mean and its $v$-derivative, where needed, have the integrable radial majorant $C\rho^{-s}$.

Away from the pole write $y=r\eta$, $a=v(\eta)$. The sign is $\sgn(a-c)$ with $c=A+(v-A)r^{-d}$. Let $k_s$ be the kernel $[1-\psi(|y-\xi|)]|y-\xi|^{-7-s}$.
Integrating first over the two block directions, differentiation in $v$ gives
\[
 \partial_v\int S_d k_s\,d\sigma
 =\int S_d\,\partial_v k_s\,d\sigma
 -2r^{-d}w(c)\E(k_s\mid a=c).
\]
The last term is zero outside $0<c<1$. The formula follows by differentiating a primitive in $a$, since $w(a)\E(k_s\mid a)$ extends continuously by zero at both endpoints. On $r\le2$ the radial majorant for this derivative is $C(r^6+r^{6-d})$, integrable even at the vertex. For $r>2$, the kernel differs from $r^{-7-s}$ by $O(r^{-8-s})$, and its $v$-derivative is $O(r^{-8-s})$ on compact interior meridians. Also $S_d$ differs from $S_\infty$ on a latitude interval of length at most $r^{-d}$. The renormalized integral at infinity and its derivative are therefore dominated, after angular integration, by $C(r^{-2-s}+r^{-1-s-d})$. All bounds are locally uniform in $A$. Dominated convergence proves the assertions; at an observation axis one uses only the estimates for the undifferentiated integral.
\end{proof}

We can now give a criterion entirely at the balanced aperture. Put
$D_d(s,A,v)=R_d(s,A,v)-R_d(s,A,A)$. On either meridian write $v=B+ez$, $z>0$, with $e=1$ towards $v=1$ and $e=-1$ towards $v=0$. Orient the functions by multiplying them by $-e$, so that their value at the observation point is positive. Denote their quadratic curvature and the gain on passing to the quartic, both relative to the link, by
\[
 h_0(z)=-eD_2(0,B,v),\qquad
 \mathcal B_0(z)=-e\,[D_2(0,B,v)-D_4(0,B,v)].
\]
Dependence on the side is understood. Their difference satisfies
\begin{equation}\label{eq:ratio}
 \frac{\mathcal B_0(z)-h_0(z)}{z}=\frac{D_4(0,B,v)}{v-B}.
\end{equation}

\begin{proposition}\label{prop:transfer}
Suppose $\mathcal B_0(z)>h_0(z)$ throughout both meridians up to the axes, and the left-hand side of \eqref{eq:ratio} has positive lower limit as $z\downarrow0$ on each side. Then Theorem~\ref{th:main} holds.
\end{proposition}
\begin{proof}
The function $D_4(s,A,v)/(v-A)$ extends continuously across $v=A$ by the integral of $\partial_vR_4(s,A,A+t(v-A))$ for $0\le t\le1$, by Lemma~\ref{lem:c1}. At $(s,A)=(0,B)$ it is positive on the compact meridian $[0,1]$. It therefore remains positive there at $(s,B_s)$ for every sufficiently small $s>0$. Indeed, a failure would give $s_j\to0^+$ and, up to a subsequence, $v_j\to v_\infty$ where the limiting quotient is nonpositive. Since $F_{4,B_s}$ has zero curvature at the stationary link, $D_4(s,B_s,v)=H_s[G_{B_s}](\xi(v))$. This gives \eqref{eq:needed}, and \eqref{eq:cal} proves the assertion.
\end{proof}

\medskip\section{Quadratic and quartic curvatures at the balanced aperture}\label{sec:critical}
We henceforth fix $B$ as in \eqref{eq:median}. On the meridian $v=B+ez$ put
$T=z/C$ for $e=1$ and $T=z/B$ for $e=-1$; thus $T$ runs from the link to the axis as it increases from zero to one.

\subsection*{Curvature of the quadratic leaves}
The ray formula \eqref{eq:ray} reduces the quadratic contribution to an integral that depends on the direction of $\nabla Q_B$.

\begin{lemma}\label{lem:four-curvature}
For $v=B+CT$ and $0<T\le1$,
\begin{equation}\label{eq:threecost}
 \frac{h_0(CT)}{CT}\le\frac{d_*}{1+2T/5}.
\end{equation}
For $v=B-BT$ and $0<T\le1$,
\begin{equation}\label{eq:curvaffine}
 H(T):=\frac{h_0(BT)}{BT}\le
 \begin{cases}
 d_*+5T/24-5T^2/48,&0<T\le1/3,\\
 3/5+T/10,&1/3\le T\le1 .
 \end{cases}
\end{equation}
\end{lemma}
\begin{proof}
\emph{Dependence on the normal.}
Let $\mathsf M_B=\operatorname{diag}(-C I_3,B I_4)$.
Subtract \eqref{eq:ray} at the link from its value at $\xi=\xi(v)$
and let $s\to0^+$. For $t=|\mathsf M_B\xi\cdot\eta|$ the difference
quotient $(t^{-s}-1)/s$ tends to $-\log t$ and is bounded by
$t^{-1/4}|\log t|$ for $s\le1/4$. This bound is integrable since
$|\mathsf M_B\xi|\ge\min(B,C)>0$ and a coordinate on the sphere has bounded
density. The factor $|Q_B(\eta)/2|^s$ tends to one almost everywhere
and is bounded by one. The length of $\mathsf M_B\xi$ contributes a constant
whose signed average vanishes by balance. Its direction has
three-block fraction
\[
 P(v)=\frac{C^2v}{B^2+(C-B)v},\qquad P(B)=C.
\]
Consequently $h_0(z)=-e[L(P(v))-L(C)]$, where
\[
 L(P)=\int_0^1\sgn(B-a)\Lambda_a(P)w(a)\,da,\qquad
 \Lambda_a(P)=\E\log\big|\sqrt{aP}\,U+\sqrt{(1-a)(1-P)}\,V\big|.
\]
The full spherical mean of $\Lambda_a$ is independent of $P$ by
rotation invariance, so the derivative of $L$ is twice the contribution
of the cap $a<B$.

To evaluate that contribution, let $\varphi(x)=\E\log|x-V|$ for
the semicircle law. Its Stieltjes transform is
$\E(z-V)^{-1}=2(z-\sqrt{z^2-1})$, with the root asymptotic to $z$.
Indeed, for $z>1$ set $V=\cos\theta$ and use
$\sin^2\theta=(1-z^2)+(z-\cos\theta)(z+\cos\theta)$ and
$\int_0^\pi(z-\cos\theta)^{-1}d\theta=\pi/\sqrt{z^2-1}$,
the latter by $t=\tan(\theta/2)$. Analytic continuation gives the
transform in the upper half-plane. Taking real parts at the boundary
gives $\varphi'(x)=2(x-\sqrt{(x^2-1)_+})$ for $x\ge0$.
For this last step, the regularizations $\E\log|x+i\epsilon-V|$
converge uniformly by the bounded density and
$\int_{\R}\log(1+\epsilon^2/t^2)\,dt=O(\epsilon)$.
Their derivatives, given by the transform, converge locally uniformly
also at the endpoint $x=1$.

Put $r^2=aP/[(1-a)(1-P)]$. Averaging first in the uniform variable
yields $\Lambda_a(P)=\tfrac12\log[(1-a)(1-P)]+\int_0^1\varphi(ru)\,du$.
Since $r'/r=1/[2P(1-P)]$,
\[
 2P(1-P)\Lambda_a'(P)=g(r)-P,\qquad
 g(r)=\frac1r\int_0^r x\varphi'(x)\,dx
     =\frac23r^2[1-(1-r^{-2})_+^{3/2}].
\]
Here $0\le g\le1$ by Bernoulli's Inequality, which justifies
differentiation under the $a$-integral on compact $P$-intervals in
$(0,1)$. Multiplication by $w(a)$ gives
$g(r)w(a)=5P[a^{3/2}-P^{-3/2}(a-1+P)_+^{3/2}]/[2(1-P)]$.
Integration from zero to $B$, using the half mass of the cap, yields
\begin{equation}\label{eq:costder}
 L'(P)=\frac{b^5-P^{-3/2}(P-C)_+^{5/2}}{(1-P)^2}
             -\frac1{2(1-P)} .
\end{equation}
The median identity gives $1/(2B)-b=d_*$.

\emph{The side $v>B$.}
Writing $y^2=(P-C)/(BP)$ gives $P=C/(1-By^2)$,
$1-P=B(1-y^2)/(1-By^2)$. Inserting these in \eqref{eq:costder} and
using $5b^3-3b^5=1$ and
$2y^5-5y^2+3=(1-y)^2(2y^3+4y^2+6y+3)$ gives
\[
 -L'(P)=bC(1-By^2)\left[\frac32+
             y^2\frac{y+1/2}{(1+y)^2}\right]\le d_*.
\]
Indeed, the last fraction is at most $1/2$ and $3B>1$.
Integrate from $C$ to $P$ and use
$(P-C)/z=BC/[BC+(C-B)z]$ and $(C-B)/B>2/5$ to obtain \eqref{eq:threecost}.

\emph{The side $v<B$.}
Here a primitive in \eqref{eq:costder} is
$b^5/(1-P)+\tfrac12\log(1-P)$. For $k=B/C$ and $u=T/(1+kT)$,
we have $1-P=B/(1-u)$ and hence $h_0(BT)=-b^3u-\tfrac12\log(1-u)$.
In the variable $T$, the logarithm is the difference of
$\log(1+kT)$ and $\log(1+(k-1)T)$. Writing this difference as an
integral and subtracting $d_*=1/(2B)-b$ gives
\[
 H(T)-d_*=\frac{T}{1+kT}\left(\frac{I(T)}{4B}-kd_*\right),
 \qquad I(T)=\int_0^1\frac{2x}{1+(k-x)T}\,dx .
\]
The integrand is convex in $T$ and $I(0)=1$, so only an upper bound at
$T=1$ is needed. Since $k>7/10$, integration with $k=7/10$ gives
\[
 I(1)\le\frac{34}{5}\operatorname{arctanh}(5/12)-2
 \le 1+\frac{163}{9072}<1+\frac1{48}.
\]
The middle bound uses $\operatorname{arctanh}y\le
y+y^3/3+y^5/[5(1-y^2)]$, obtained from its power series.
It follows by the chord inequality that $I(T)\le1+T/48$.
Moreover $1/(4B)-kd_*=1/(4b^2)-3b^3/2<5/24$, by monotonicity
and $b>9/14$. As $B>2/5$, the bracket in the formula for $H(T)-d_*$
is at most $(5/24)(1+T/16)$. Now
$(1+T/16)(1+3T/5)\le1+7T/10$ on $[0,1]$, which gives $H(T)\le d_*+5T/[24(1+3T/5)]$. For the affine bound in \eqref{eq:curvaffine} use
$d_*<17/30$ and $20-53T+36T^2>0$, whose discriminant is negative.
For $T\le1/3$, $(1-T/2)(1+3T/5)\ge1$ gives the small-$T$ bound in the same formula.
The inequalities extend to the observation axes by continuity of $D_2$ in Lemma~\ref{lem:c1}.
\end{proof}

\subsection*{Reflection between reciprocal spheres}
The gain $\mathcal B_0$ in \eqref{eq:ratio} is obtained by pairing reciprocal radii. For source latitude $a$ set
\[
 K_r(v,a)=\E\big(1+r^2-2r[\sqrt{va}\,U+
                 \sqrt{(1-v)(1-a)}\,V]\big)^{-7/2},\qquad
 f_r(v,a)=w(a)K_r(v,a).
\]
For $0<r<1$ extend $f_r$ by zero outside $[0,1]$. This extension is
continuous since $w$ vanishes at the endpoints. The angular kernel $K_r$ itself is defined at the axes by the same expectation.

Fix $0<r<1$. The outer sphere of radius $1/r$ and the inner sphere of radius $r$ determine the two intervals we shall compare. Write $l=|y|$ and orient $Q_B,G_B$ by multiplying both by $-e$. At $\xi$ both oriented functions have value $z$, so their respective curvature signs are $S_d(y)=\sgn(z-e l^d(a-B))$ for $d=2,4$. For $l>1$ the difference $S_2-S_4$ equals $2$ on $zl^{-4}<e(a-B)<zl^{-2}$. For $l<1$ it equals $-2$ on $zl^{-2}<e(a-B)<zl^{-4}$. The difference vanishes elsewhere up to the null level sets. Thus the coordinate $e(a-B)/z$ ranges over $[r^4,r^2]$ on the outer sphere and $[r^{-2},r^{-4}]$ on the inner sphere. Pair the two intervals in opposite order by setting
\begin{equation}\label{eq:paircoords}
 q'=\frac{r^2+r^4-q}{r^6},\qquad a=B+ezq,\qquad a'=B+ezq',
 \qquad r^4<q<r^2 .
\end{equation}

For $\epsilon>0$, denote by $\mathcal B_{0,\epsilon}(z)$ the integral of $S_2-S_4$ with the regularized kernel $\omega_6^{-1}(|\xi-y|^2+\epsilon^2|y|)^{-7/2}$. Let $f_{r,\epsilon}$ be $f_r$ with $\epsilon^2r$ added to the squared denominator in $K_r$, with the same zero extension. Under $l=1/r$ the squared denominator is $r^{-2}$ times its value at $l=r$. The outer radial weight becomes $dr/r$, and the inner weight is $r^6\,dr$. Since $|dq'/dq|=r^{-6}$, we obtain
\begin{equation}\label{eq:smoothpair}
 \frac{\mathcal B_{0,\epsilon}(z)}z
 =2\int_0^1\frac{dr}{r}\int_{r^4}^{r^2}
 \big[f_{r,\epsilon}(v,a)-r f_{r,\epsilon}(v,a')\big]\,dq .
\end{equation}
To remove the regularization we first compare the two angular densities at each parameter of the Laplace integral.

\begin{lemma}\label{lem:angular-reflection}
For the latitudes $(a,a')$ in \eqref{eq:paircoords} and every $\beta>0$, put
\[
 p_\beta(u)=w(u)\E\exp\big(\beta[\sqrt{vu}\,U+\sqrt{(1-v)(1-u)}\,V]\big)
\]
on $[0,1]$, and put $p_\beta=0$ outside. If $e=1$, then $p_\beta(a')\le p_\beta(a)$.
If $e=-1$, then $r p_\beta(a')\le p_\beta(a)$, and for $a'>0$ also
$p_\beta(a')\le p_\beta(a)[(1-a')/(1-a)]^{1/3}$.
\end{lemma}
\begin{proof}
We work with interior latitudes. An incoming latitude outside $(0,1)$ has zero density, and we treat observation axes by continuity. At fixed $\beta>0$ the density to be compared at $u=a,a'$ is
$w(u)\E e^{AU+DV}$, where $A=\beta\sqrt{vu}$ and
$D=\beta\sqrt{(1-v)(1-u)}$. Put $\theta(u)=\arcsin\sqrt u$.
Up to a positive factor independent of $u$, this density is
\[
 e^{\beta\cos(\theta(u)-\theta(v))}P(u),\qquad
 P(u)=(1-e^{-2A})De^{-D}I_1(D).
\]
Here $I_0(x)=\pi^{-1}\int_0^\pi e^{x\cos t}dt$ and $I_1=I_0'$.
Integration by parts gives $I_1'=I_0-I_1/x$ and
$\E e^{DV}=2I_1(D)/D$, while $\E e^{AU}=\sinh(A)/A$.

We first show that the reflected point is farther from the observation latitude in the meridian angle.
In fact $k_0=(q'-1)/(1-q)\ge r^{-2}$ since
$r^6(1-q)(k_0-r^{-2})=(r^2-q)(1-r^4)$. For $e=-1$ we have
$a'<v<a\le B<1/2$, where $\theta'$ is decreasing, and
$v-a'=k_0(a-v)$. For $e=1$ write $d=v-a$, $a'=v+k_0d$.
For $0\le h\le1$, squaring
$k_0\theta'(v+k_0dh)\ge\theta'(v-dh)$ reduces it to
\[
 (1-1/k_0)v(1-v)\ge (1-2v)dh.
\]
For $v\ge1/2$ the right-hand side is nonpositive. For $B\le v<1/2$
use $d\le2(1-r^2)(v-B)$, $1-1/k_0\ge1-r^2$ and
$2(1-2v)\le1-v$, since $v\ge1/3$. Integration in $h$ proves the
angular comparison. It remains to control $P$, for which
\begin{equation}\label{eq:trfactor}
 (\log P)'(u)=\frac{A}{u(e^{2A}-1)}
 -\frac{D}{2(1-u)}\left(\frac{I_0(D)}{I_1(D)}-1\right).
\end{equation}

First let $e=1$. The ratio $R=I_1/I_0$ solves
$R'=1-R^2-R/x$. A supersolution with the same zero limit at the
origin is $Y=\tanh(x/2)$, since $\sinh x\ge x$.
The integrating factor $x\exp\int_0^x(R+Y)$ for the difference gives
$R\le Y$. It follows that
\[
 (\log P)'(u)\le\frac{A}{u(e^{2A}-1)}
                  -\frac{D}{(1-u)(e^D-1)}\le0
 \qquad(v,u\ge B).
\]
Here $2A\ge D$, $1-u\le2u$, and $x/(e^x-1)$ decreases.
Consequently the outer density dominates its reflected partner.

Now let $e=-1$. The bound we need for the reciprocal ratio is
\[
 \frac{I_0(x)}{I_1(x)}\le \coth(kx)+\frac{g}{x},\qquad
 g=\frac12+\frac{2k}{9},\qquad 0<k<\frac34.
\]
To verify it, the reciprocal ratio satisfies $R'=1-R^2+R/x$.
For the proposed upper bound $F$, with $y=kx$, its residual is
\[
 x^2(F'-1+F^2-F/x)=
 \frac{1-k}{k^2}\left(\frac{y}{\sinh y}\right)^2+
 \frac49 y\coth y-g(2-g).
\]
The first coefficient is at least $4/9$, and $g(2-g)\le8/9$.
The inequality $(y/\sinh y)^2+y\coth y\ge2$ follows from the
decrease of $L(y)/y$, where $L(y)=(\log\E e^{yU})'=\coth y-1/y$.
Indeed integration by parts in the uniform law gives
$L(y)/y=(1-\E_y U^2)/2$ under the exponential tilt. The law of
$|U|$ is weighted by $\cosh(y|U|)$, so $\partial_y\E_y U^2$
is the covariance of the increasing functions $u^2$ and $u\tanh(yu)$.
Thus $F$ is a supersolution, and it dominates near zero because
$I_1'(0)=1/2$ gives $xR(x)\to2<1/k+g$. Scalar comparison proves the bound.

Take $k=\sqrt{Bu/[C(1-u)]}$. Then $A/D\le k$,
$u/(1-u)\le k$ and $k\le B/C<3/4$. Using the decrease of
$x/(e^{2x}-1)$ gives
$(1-u)A/[u(e^{2A}-1)]\ge D/(e^{2kD}-1)$.
This cancels the term $D(\coth(kD)-1)/2$ in \eqref{eq:trfactor}.
Hence $-(1-u)(\log P)'\le g/2\le1/3$, which upon integration
from $a'$ to $a$ gives the asserted cube-root bound for $p_\beta$, since the exponential factor also decreases under reflection.
We also need to absorb the prefactor into the radial weight $r$.
Set $p=kC/B\le1$. Since $(B-u)/(1-u)=B(1-p^2)$ and $k\le3p/4$,
\[
 -8(B-u)(\log P)'\le 2B(1-p^2)(1+p/3)
 \le 2B(1+1/36)<1,
\]
where $p/3-p^2\le1/36$. Integrating once more gives
$P(a')/P(a)\le(q'/q)^{1/8}\le r^{-1}$. This proves the
nonnegativity of the pair at each Laplace parameter. The bounds extend to the observation axes by continuity of the angular averages for each $r<1$.

\end{proof}

\begin{proposition}\label{prop:quarticmirror}
With the latitudes in \eqref{eq:paircoords},
\begin{equation}\label{eq:pairbridge}
 \frac{\mathcal B_0(z)}z
 =2\int_0^1\frac{dr}{r}\int_{r^4}^{r^2}
          [f_r(v,a)-r f_r(v,a')]\,dq .
\end{equation}
The paired density is nonnegative on both sides. For $e=1$ one has
$f_r(v,a')\le f_r(v,a)$, hence
\begin{equation}\label{eq:mirrorlower}
 \frac{\mathcal B_0(z)}z\ge 2\int_0^1
 \frac{1-r+r\mathbf1_{\{r^2<z/C\}}}{r}
 \int_{r^4}^{r^2}f_r(v,B+zq)\,dq\,dr .
\end{equation}
For $e=-1$ one has, whenever $a'>0$,
\begin{equation}\label{eq:trcube}
 f_r(v,a')\le f_r(v,a)\left(\frac{1-a'}{1-a}\right)^{1/3}.
\end{equation}
\end{proposition}
\begin{proof}
Using $h^{-7/2}=\Gamma(7/2)^{-1}\int_0^\infty t^{5/2}e^{-th}\,dt$, we can write the angular difference in \eqref{eq:smoothpair} as
\[
 \frac{1}{\Gamma(7/2)}\int_0^\infty
 t^{5/2}e^{-t(1+r^2)}e^{-\epsilon^2rt}
 \big[p_{2rt}(a)-r p_{2rt}(a')\big]\,dt .
\]
The bracket is nonnegative by Lemma~\ref{lem:angular-reflection}. Monotone convergence as $\epsilon\downarrow0$ therefore gives the right-hand side of \eqref{eq:pairbridge}. On the curvature
side it suffices to compare with the centered smoothing
$(|\xi-y|^2+\epsilon^2)^{-7/2}$. For $\rho=|\xi-y|<1/2$ their
difference in $L^1$ is at most
$C\epsilon^2\int_0^{1/2}\rho^7(\rho^2+\epsilon^2)^{-9/2}\,d\rho
=O(\epsilon)$. The centered smoothing gives the usual principal value. On each small sphere the tangent half-space has zero signed mean, and
the angular error of the smooth leaf is $O(\rho)$ since its gradient is
nonzero. Off the pole the kernel is locally bounded; for $l>2$ the
difference of the two curvature signs is supported in a latitude interval of
length at most $zl^{-2}$, giving an integrable $O(l^{-3})$ radial tail.
The principal value obtained in this way is $\mathcal B_0(z)$, since the conical tails cancel and the two renormalized curvatures agree at the link. This proves \eqref{eq:pairbridge}. Integrating the inequalities of Lemma~\ref{lem:angular-reflection} in the Laplace parameter gives the stated angular bounds on both sides. For $e=1$ and $r^2<z/C$, the incoming density is zero since $q'\ge r^{-2}$, giving the additional term in \eqref{eq:mirrorlower}.
\end{proof}

\medskip\section{Preliminary estimates for the angular kernel}\label{sec:prelim}
The formula at reciprocal radii reduces the remaining comparison to estimates for $K_r(v,a)$. Away from the link we need a lower bound that retains the difference of the two latitudes. Near the link we compare with $K_r(B,B)$ and integrate the error against the radial measure of the limiting pair.

\begin{lemma}\label{lem:orbitalbase}
Let $0<r<1$ and $0\le v,a\le1$. Set $D=1+r^2$,
$t=2r\sqrt{va}/D$, $y=2r\sqrt{(1-v)(1-a)}/D$, and
\[
 \mathcal W=(1-r^2)^4+8r^2(1-r^2)^2(v+a-2va)+16r^4(v-a)^2 .
\]
Then
\begin{equation}\label{eq:orbitaltwo}
 \frac{\mathcal W K_r(v,a)}{\sqrt D}\ge
 1+\frac58t^2-\frac{y^2+3y^4}{32}.
\end{equation}
With $d=\sqrt{va}+\sqrt{(1-v)(1-a)}$ one also has
\begin{equation}\label{eq:poissonorbit}
 (1-r^2)K_r(v,a)\ge
 \left(1+\frac{4r^2(1-d^2)}{(1-r^2)^2}\right)^{-7/4}.
\end{equation}
\end{lemma}
\begin{proof}
Write $a_1=Dt$, $a_2=Dy$. Integrating first in the uniform coordinate gives
\begin{equation}\label{eq:jdivdiff}
 K_r(v,a)=\frac{F(D-a_1)-F(D+a_1)}{5a_1},\qquad
 F(x)=\E(x-a_2V)^{-5/2}.
\end{equation}
For $x>a_2$ put $\lambda=a_2/x$. The substitution
$V=(u+\lambda)/(1+\lambda u)$ multiplies the semicircle density by
$(1-\lambda^2)^{3/2}/(1+\lambda u)^3$ and gives
\[
 F(x)=\frac{A(a_2/x)}{\sqrt{x}(x^2-a_2^2)},\qquad
 A(\lambda)=\E(1+\lambda V)^{-1/2}=\sum_{j\ge0} c_j\lambda^{2j}.
\]
By symmetry and the binomial series, $c_j>0$, $c_0=1$ and
$c_1=(3/8)\E V^2=3/32$. The product of the two quadratic denominators
in \eqref{eq:jdivdiff} is $\mathcal W$.
For $p=1/2+2j$ set $d_p=((1-t)^{-p}-(1+t)^{-p})/t$ and
$s_p=(1-t)^{-p}+(1+t)^{-p}$. The contribution of $c_j$ to
$\mathcal W K_r/\sqrt D$ is
\[
 \frac{c_jy^{2j}}5\left[(1+t^2-y^2)d_p+2s_p\right].
\]
All terms are nonnegative and the series converge since $t+y\le2r/D<1$.
Keep only $j=0,1$. For $j=0$ the bounds $d_{1/2}\ge1+5t^2/8$ and
$2s_{1/2}\ge4+3t^2/2$ give a contribution of at least
$1+5t^2/8-y^2/5-t^2y^2/8$. For $j=1$ use $d_{5/2}\ge5$ and
$2s_{5/2}\ge4+2t^2$, again from the binomial series. This contributes
at least $3y^2(9-5y^2+7t^2)/160$. Adding proves
\eqref{eq:orbitaltwo}; the divided differences at $t=0$ have their
continuous values.

To prove \eqref{eq:poissonorbit} we use the Poisson kernel. On the unit ball in $\R^7$ its restriction to functions
of latitude has density $P_r(v,a)=(1-r^2)K_r(v,a)$ relative to
$w(a)\,da$ \cite[Theorem~1.17]{ABR}. It is symmetric, preserves
constants, and $P_tP_\rho=P_{t\rho}$ since $x\mapsto u(\rho x)$ is
the harmonic extension of $\xi\mapsto u(\rho\xi)$. By continuity of the
kernels at radii less than one this identity holds at all latitudes. Thus
$P_r(u,u)=\int P_{\sqrt r}(u,a)^2w(a)\,da\ge1$ by Cauchy--Schwarz.
For $d>0$ write $u=\sqrt{va}/d$ and factor the angular denominator as
$h[1+t^2-2t(uU+(1-u)V)]$ with $ht=rd$, $h(1+t^2)=1+r^2$.
For $L=[1+4r^2(1-d^2)/(1-r^2)^2]^{1/2}$ the larger root is
$h=[1+r^2+(1-r^2)L]/2\in[1,L]$ and $t=rd/h<1$.
Since $h(1-t^2)=(1-r^2)L$, the diagonal estimate gives
$(1-r^2)K_r(v,a)\ge h^{-5/2}/L\ge L^{-7/2}$. Continuity covers $d=0$.
\end{proof}

\begin{lemma}\label{lem:link-gaussian}
Put $k_*=BU+CV$ and $d=\sqrt{va}+\sqrt{(1-v)(1-a)}$. There is a constant
$\eta\in(0,1)$ with $\eta/(1+b)>3/5$ such that
\[
 K_r(v,a)\ge\eta\,\E\big(1+r^2-2rk_*+2r(1-d)\big)^{-7/2}
\]
for all $0<r<1$ and $v,a\in[0,1]$.
\end{lemma}
\begin{proof}
Let $U_n$ be the first coordinate on the unit sphere in $\R^n$ with
normalized surface measure, $F_n(t)=\E e^{t(U_n-1)}$, and
$R_n(t)=(\log\E e^{tU_n})'$. The exponential tilt by $e^{tU_n}$ gives
$R_n'=\operatorname{Var}_t U_n>0$. Integrating the derivative of
$e^{tu}(1-u^2)^{(n-1)/2}$ on $[-1,1]$ yields
\[
 R_n/t=(1-\E_t U_n^2)/(n-1),\qquad
 R_n'=1-R_n^2-(n-1)R_n/t .
\]
Also $R_n/t$ decreases. Indeed the law of $|U_n|$ is tilted by
$\cosh(tu)$, and the derivative of its second moment is the covariance
of the increasing functions $u^2$ and $u\tanh(tu)$.

We locate the minimum of $g(u)=\log F_3(tu)+\log F_4(t(1-u))$.
Its derivative is $t[R_3(tu)-R_4(t(1-u))]$ and is increasing.
For $P_n(x)=R_n(nx)$ and $Q_n(x)=R_n((n-1)x)$ the two rescaled equations are
\[
 P_n'=n(1-P_n^2)-(n-1)P_n/x,\qquad
 Q_n'=(n-1)(1-Q_n^2-Q_n/x).
\]
In the first equation $P_4$ is a subsolution for $n=3$ since
$R_4(4x)/(4x)-R_4'(4x)\ge0$. In the second $Q_3$ is a subsolution for
$n=4$ since $R_3'(2x)>0$. For the differences $P_3-P_4$ and $Q_4-Q_3$
the integrating factors are respectively
$x^2\exp(3\int_0^x(P_3+P_4))$ and $x^3\exp(3\int_0^x(Q_3+Q_4))$.
The products vanish at zero, giving $P_3\ge P_4$ and $Q_4\ge Q_3$.
It follows that $g'(2/5)\le0\le g'(3/7)$, so the minimum $u_*$ and $B$
both lie in $[2/5,3/7]$.

A single variance bound applies to both spherical factors. Under the tilt,
$h=z(1-U_n)$ has density proportional to
$h^{p-1}e^{-h}(2-h/z)^{p-1}\one_{0<h<2z}$, $p=(n-1)/2$.
For $n=3,4$ the factor multiplying the gamma density is nonincreasing
on its support, with value zero beyond it. Weighting by this factor
can only lower the second moment. Thus $z^2R_n'(z)=\operatorname{Var}_z h
\le p(p+1)$. For $2/5\le u\le3/7$ this gives
\[
 g''(u)\le\frac2{u^2}+\frac{15}{4(1-u)^2}
 \le\frac{25}{2}+\frac{735}{64}<24 .
\]
Taylor's Formula yields $g(B)-g(u_*)\le12(3/7-2/5)^2<1/100$.
Put $u=\sqrt{va}/d$ when $d>0$. Since $d\le1$ and $F_n$ is decreasing,
\[
 F_3(t\sqrt{va})F_4(t\sqrt{(1-v)(1-a)})
 \ge e^{g(u)}\ge e^{-1/100}F_3(tB)F_4(tC).
\]
The same holds for $d=0$ by continuity. Apply this bound at $t=2r\gamma$ for each integration parameter $\gamma>0$ in the Laplace representation of the power $-7/2$. Restoring the factor
$e^{td}$ increases the link denominator by $2r(1-d)$. This gives the
claim with $\eta=e^{-1/100}>99/100$, since $b<13/20$.
\end{proof}

As $z\downarrow0$ the two latitudes $a,a'$ of a fixed pair approach $B$. In the quotient \eqref{eq:pairbridge} the pair then has density $(1-r)w(B)K_r(B,B)$. The length of the $q$-interval is $r^2(1-r^2)$. These factors give the radial measure used below.

\begin{lemma}\label{lem:linkmoments}
The link density $H_*(r)=2w(B)r(1-r)(1-r^2)K_r(B,B)$ satisfies
\[
 J:=\int_0^1H_*(r)\,dr=(1+b)^{-1},\qquad
 H_*(r)\le\min\{1,3r/(1+r)\}.
\]
Let $r$ have probability density $H_*/J$. Given $r$, tilt the law of $(U,V)$ by
$(1+r^2-2rk_*)^{-7/2}/K_r(B,B)$ and take $q$ uniform on $(r^4,r^2)$,
independently of the angles. With $\tau=(1-r)/\sqrt r$, $\ell=2(1-k_*)$ and
$M=r\,\E[((1-q)/(1-r))^2\mid r]$, this law satisfies
\[
 J\E q\le4/15,\qquad J\E\!\left[M\frac{\tau^2}{\tau^2+\ell}\right]\le6/7.
\]
\end{lemma}
\begin{proof}
For $k=k_*<1$ the denominator $1+r^2-2rk$ is $r(\tau^2+\ell)$.
Using $d\tau/dr=-(1+r)/(2r^{3/2})$ gives
\[
 \int_0^1\frac{r(1-r)(1-r^2)\,dr}{(1+r^2-2rk)^{7/2}}
 =2\int_0^\infty\frac{\tau^2\,d\tau}{(\tau^2+\ell)^{7/2}}
 =\frac1{15(1-k)^2},
\]
where the last substitution is $\tau=\sqrt\ell\tan\theta$.
Averaging $(1-k_*)^{-2}$ first in $U$ gives $[(1-CV)^2-B^2]^{-1}$.
Partial fractions and the semicircle transform
$\E(x-\alpha V)^{-1}=2/(x+\sqrt{x^2-\alpha^2})$ at $x=1\pm B$,
$\alpha=C$, give
\[
 J=\frac{2w(B)}{15B}\left[\frac1C-\frac1{(1+b)^2}\right]=\frac1{1+b}.
\]
The integrands in this calculation are nonnegative; the transform at $x=\alpha$ follows by monotone convergence.

To bound $H_*$ from above, take $x=1+r^2-2Br$, $\alpha=2Cr$, $D=x+\alpha$,
$\Delta=\sqrt{x^2-\alpha^2}=(1-r)\sqrt D$. In \eqref{eq:jdivdiff}
retain only the first term. Jensen's Inequality and the semicircle transform give
\[
 F(x)=\frac{\E(1+(\alpha/x)V)^{-1/2}}{\sqrt x\,\Delta^2}
 \le\frac{\sqrt{2/(x+\Delta)}}{\Delta^2}
 =\frac{2}{(1-r)^2D(1-r+\sqrt D)} .
\]
It follows that $H_*\le [2w(B)/(5B)](1+r)/[D(1-r+\sqrt D)]\le1$.
Indeed $D=C(1+r)^2+B(1-r)^2$ gives $1-r+\sqrt D\ge2\sqrt C$, and
\[
 \frac D{1+r}=1+r-4B+\frac{4B}{1+r}
 \ge4b(1-b)>\frac9{10}>\frac{w(B)}{5B\sqrt C},
\]
where $1/2<b<13/20$ and $B/C>7/10$.
For the second bound write $1+r^2-2rk_*=BD_3+CD_4$ with
$D_n=1+r^2-2rU_n$ and independent coordinates. The arithmetic--geometric
mean inequality with weights $3/7,4/7$ and the Poisson integral of the
constant in the respective dimensions give
\[
 K_r(B,B)\le (3/(7B))^{3/2}(4/(7C))^2\E D_3^{-3/2}\E D_4^{-2}
 =(3/(7B))^{3/2}(4/(7C))^2(1-r^2)^{-2}.
\]
After multiplication by $2w(B)$ the prefactor is
$15(3/7)^{3/2}(4/7)^2/(2BC)<3$ since $BC>6/25$.

Finally $J\E q\le\int_0^1(r^2+r^4)\,dr/2=4/15$. Write
$q=(1-\lambda)r^2+\lambda r^4$ with $\lambda$ uniform on $(0,1)$.
Then $M=r(1+r)^2A_0(r)$, $A_0=1+r^2+r^4/3$. In $(\tau,U,V)$ the
unnormalized measure has density $4w(B)\tau^2(\tau^2+\ell)^{-7/2}$ relative
to $d\tau$ and the original angular law. Integrating by parts in $\tau$
at fixed $\ell>0$ gives
\[
 J\E\!\left[M\frac{\tau^2}{\tau^2+\ell}\right]
 =\frac37\int_0^1 H_*(r)R(r)\,dr,\qquad
 R=M-\frac{2r(1-r)}{3(1+r)}M'.
\]
The boundary term $\tau^3M(\tau^2+\ell)^{-7/2}$ vanishes at both ends.
Moreover $3\tau^2 R=\partial_\tau(\tau^3M)\ge0$, since
$\tau^3M=\tau\E_\lambda(1-q)^2$ and $q$ decreases with $\tau$ at fixed
$\lambda$. Thus Tonelli applies and we may use $H_*\le1$. A last
integration by parts in $r$, with $2r(1-r)/(1+r)$ zero at both ends, gives
\[
 3\int_0^1 R\,dr=\int_0^1 r[(1+r)^2+4]A_0(r)\,dr
 \le\int_0^1(3+x)(1+x+x^2/3)\,dx=23/4<6,
\]
using $(1+r)^2\le2(1+r^2)$ and $x=r^2$. This proves the second moment bound.
\end{proof}

\medskip\section{The sign on the two sides of the link}\label{sec:meridians}
On the side $v>B$ reflection decreases the angular density. On $v<B$ the angular density may increase, and it matters which radii have an incoming point. We treat the two sides in turn.
\subsection*{The three-dimensional axis}
\begin{proposition}\label{prop:plus}
For $v=B+z$ one has $\mathcal B_0(z)>h_0(z)$ throughout $0<z\le C$, and
$\liminf_{z\downarrow0}(\mathcal B_0(z)-h_0(z))/z>0$.
\end{proposition}
\begin{proof}
Put $T=z/C$. The estimate \eqref{eq:threecost} gives
$h_0(z)/z\le d_* /(1+2T/5)$. We first compare with the link near
$T=0$, and then retain the additional flux from the radii at which
the reflected point has passed the axis.

For $T\le1/5$ both $a=B+zq$ and $v$ belong to $[B,C]$. Writing
$\theta(u)=\arcsin\sqrt u$ and $d=\cos(\theta(v)-\theta(a))$, we have
\[
 2(1-d)\le (\theta(v)-\theta(a))^2
 \le Z(1-q)^2,\qquad Z=\frac{C}{4B}T^2<\frac38T^2 .
\]
Also $w(a)/w(B)\ge 1-Tq/2$. Indeed $(1-y)\sqrt{1+Cy/B}$ is concave,
equals one at zero and exceeds $9/10$ at $y=1/5$. In the link comparison
of Lemma~\ref{lem:link-gaussian}, the relative increase in the squared
distance is bounded by
\[
 \delta=Zr\left(\frac{1-q}{1-r}\right)^2\frac{\tau^2}{\tau^2+\ell}.
\]
Keeping the term with weight $1-r$ in \eqref{eq:mirrorlower} and
using $(1-u)(1+\delta)^{-7/2}\ge1-u-7\delta/2$ for $0\le u\le1$,
we obtain from Lemma~\ref{lem:linkmoments}
\[
 \frac{\mathcal B_0(z)}z
 \ge \eta J-\tfrac12 TJ\E q-\tfrac72 J\E\delta
 >\frac35-\frac{2T}{15}-\frac98 T^2 .
\]
Here $J\E\delta\le6Z/7$. The difference between the last expression
and $d_* /(1+2T/5)$ is concave. It is positive at zero since $d_*<17/30$,
leaving a positive lower limit at the link. At $1/5$ use $(27/25)(3/5-2/75-9/200)>17/30$.

Suppose now $1/5\le T\le1$. Write $x=r^2$ and use a uniform
$\lambda\in[0,1]$ to parametrize the outgoing interval by
$q=x-x(1-x)\lambda$. In \eqref{eq:orbitaltwo}, put $u=2\sqrt x/(1+x)$,
so $t^2=u^2va$ and $y^2=u^2(1-v)(1-a)$. Since $v\ge B+C/5>13/25$ and
$a\ge B>2/5$, we have $va>26/125$ and $(1-v)(1-a)<3/10$.
It follows that $\mathcal W K_r(v,a)\ge\sqrt{1+x}(1+u^2/9)$,
as $(5/8)(26/125)-(3/10+27/100)/32=359/3200>1/9$. The ratio
$(1+u^2/9)/\sqrt{1+u}$ decreases on $[0,1]$, as its derivative has the
sign of $3u^2+4u-9$. Thus $w(B)=(5/2)d_*$ gives
\[
 w(B)K_r(v,a)\ge\frac{25d_*}{9\sqrt2}\frac{1+\sqrt x}{\mathcal W}.
\]
We use $w(a)/w(B)\ge1-Tq$ and average the denominator with this weight.
Its mean weight is $W_T(x)=1-T(x+x^2)/2$. After normalization
$\lambda$ has density proportional to $1+k\lambda$, with
$k=Tx(1-x)/(1-Tx)\le1$. For this density
$\E(1+x\lambda)^2\le1+3x/2$. In fact
$2\E\lambda+\E\lambda^2$ increases with $k$ and equals $3/2$ at $k=1$, as one sees by integrating against $(1+k\lambda)/(1+k/2)$.
Now
\[
 \frac{\mathcal W}{(1-x)^2}=(1-x)^2+8x(v+a-2va)
                         +16x^2z^2(1+x\lambda)^2 .
\]
In the middle term use $q\ge x^2$ and $1-2v<0$. The weighted mean is
then bounded above by
\begin{align*}
 &(1-x)^2+16BCx+8(1-2B)zx(1+x^2)+16z^2x^2+8z^2x^3\\
 &\hspace{12mm}\le D_T(x):=(1+x)^2+\tfrac67Tx(1+x^2)+3T^2x^2(2+x),
\end{align*}
using $8(1-2B)C<6/7$ and $8C^2<3$ from $B>(9/14)^2$.
Jensen's Inequality for the reciprocal in \eqref{eq:mirrorlower}
therefore gives
\begin{equation}\label{eq:plusradial}
 \frac{\mathcal B_0(z)}{zd_*}\ge \frac{25}{9\sqrt2}(I_1+I_2),\qquad
 I_1=\int_0^1\frac{W_T(x)}{D_T(x)}\,dx,\quad
 I_2=\int_0^T\frac{\sqrt x+x}{D_T(x)}\,dx .
\end{equation}
Indeed $2\,dr/r=dx/x$ and the $q$-interval has length $x(1-x)$.
When $x<T$ the full outgoing flux is retained, adding the fraction
$\sqrt x/(1-\sqrt x)$ of the first integrand. Here $W_T\ge1-x$.

The one-variable estimate in Lemma~\ref{lem:radialplus} gives
$(1+2T/5)(I_1+I_2)\ge51/100$. Substitution in \eqref{eq:plusradial} completes the comparison, since $(25/(9\sqrt2))(51/100)=17/(12\sqrt2)>1$.
\end{proof}

\subsection*{The four-dimensional side near the link}
\label{sec:minusnear}
On this side $z=BT$, $v=B-z$ and $a=B-zq$.
For $x=r^2$ write
\[
 q=x-x(1-x)\lambda,\qquad
 q'=x^{-1}+(1-x)x^{-2}\lambda,\qquad a'=B-zq'.
\]
Brackets denote averaging in $0<\lambda<1$.
The reflected point is present precisely when $r>\rho=\rho_T(\lambda)$,
where
\begin{equation}\label{mn:radius}
 \rho^4=T[(1-\lambda)\rho^2+\lambda],\qquad
 t^2\le\rho\le t,\qquad t=T^{1/4}.
\end{equation}
Indeed the positive root in $\rho^2$ lies between $T$ and $\sqrt T$.
The contribution lost to reflection on $r>\rho$ and the additional
outgoing contribution on $r<\rho$ are governed by the same moment
$\langle\rho^3\rangle$.

\begin{lemma}\label{prop:otherfifth}
For $0<T\le1/3$ on $v=B-BT$, $\mathcal B_0(z)>h_0(z)$.
The lower limit of $(\mathcal B_0(z)-h_0(z))/z$ at the link is positive.
\end{lemma}
\begin{proof}
Put $k=B/(3C)<5/21$. The cube-root bound \eqref{eq:trcube}
gives $f_r(v,a')\le f_r(v,a)[1+kT(q'-q)]$ for $r>\rho$.
As the paired density is nonnegative, it is bounded below by
\[
 f_r(v,a)\big[(1-r)(1-L)+r\one_{r<\rho}\big],\qquad
 L=\min\{1,kTr(q'-q)/(1-r)\}\one_{r>\rho}.
\]
Its radial measure in \eqref{eq:pairbridge} is $2r(1-r^2)\,dr$,
since $q$ ranges over an interval of length $r^2(1-r^2)$.

We first estimate the term with $1-r$. Let
$d=\sqrt{va}+\sqrt{(1-v)(1-a)}$. Since $2B/3\le v\le a\le B<1/2$,
the density bound \eqref{eq:weight} gives $w(a)\ge w(B)$, and
\[
 2(1-d)\le(\arcsin\sqrt a-\arcsin\sqrt v)^2
 \le Z(1-q)^2,\qquad
 Z=\frac{BT^2}{4(1-T)(C+BT)}\le\frac{13}{60}T^2 .
\]
The derivative of the meridian angle is largest at $v$ in this range.
The coefficient of $T^2$ in $Z$ increases in $B$ and $T$,
and at $B=5/12$, $T=1/3$ it equals $45/208<13/60$.
In the link kernel of Lemma~\ref{lem:link-gaussian} the relative increase
of the denominator is at most
$\delta=Z(1-q)^2/(\tau^2+\ell)$, with the notation of
Lemma~\ref{lem:linkmoments}. Under the link law,
$J\E\delta\le6Z/7$, since $\tau^2=(1-r)^2/r$.
Using $(1-L)(1+\delta)^{-7/2}\ge1-L-7\delta/2$, we obtain
\begin{equation}\label{mn:balance}
 \mathcal B_0(z)/z\ge \eta J-3Z+G-D,
\end{equation}
where
\[
 G=2\Big\langle\int_0^\rho r^2(1-r^2)f_r(v,a)\,dr\Big\rangle,
 \qquad
 D=kT\Big\langle\int_\rho^1H_*(r)\frac{r(q'-q)}{1-r}\,dr\Big\rangle .
\]

\emph{The two integrals at the cutoff.} To estimate $G-D$, observe that the singular term in the majorant for $D$ is proportional to $T\lambda\int_\rho^1 r^{-2}\,dr$. The equation for the cutoff gives $T\lambda/\rho=\rho^3-T(1-\lambda)\rho$, and the term $\rho^3$ also measures the radial mass in $G$. Set $m_3=\langle\rho^3\rangle$.

On the outgoing interval $r<\rho$ we have $1-d^2\le Z(1-q)^2$ and
$q\ge r^4$. The Poisson bound \eqref{eq:poissonorbit} and Bernoulli's
inequality yield
\[
 (1-r^2)K_r(v,a)\ge[1+4Zr^2(1+r^2)^2]^{-7/4}
             \ge1-7Zr^2(1+r^2)^2 .
\]
Averaging the last error with density $3r^2/\rho^3$ gives a factor
\[
 \tfrac35\rho^2+\tfrac67\rho^4+\tfrac13\rho^6
 \le \tfrac9{25}+\tfrac27+\tfrac1{15}<\tfrac57 ,
\]
because $\rho^2\le\sqrt T<3/5$ and $\rho^4\le T\le1/3$.
Consequently, using
$2w(B)/3>15/16$, we have
\[
 G\ge\tfrac{15}{16}(1-5Z)m_3
    \ge(\tfrac{15}{16}-\tfrac{65}{64}T^2)m_3 .
\]

For the loss $D$ use $H_*\le3r/(1+r)$. After division by $3kT$, the integrand in this upper bound is
$A_\lambda(r)=1+r^2+\lambda(r^{-2}+r^4)$.
On $r\ge m=\max(t,1/2)$ all reflected points are present.
The second bound $H_*\le1$ reduces this majorant by the fraction
$(2r-1)/(3r)$. Since $\langle A_\lambda\rangle\ge3r$ by
arithmetic--geometric mean, the saving is at least
$m(1-m)\ge t(1-t)$.
A primitive of $A_\lambda$ is $r+r^3/3+\lambda(-1/r+r^5/5)$,
whose mean at $1$ is $14/15$.
At the lower endpoint \eqref{mn:radius} gives
$T\lambda/\rho=\rho^3-T(1-\lambda)\rho$. Dropping the term
$-T\lambda\rho^5/5$ bounds $D/(3k)$ by the quantity in braces below,
which is therefore nonnegative because $D\ge0$. Since $3k<5/7$, this gives
\[
 D\le\tfrac57\big\{m_3+
   T[\tfrac{14}{15}-t(1-t)-\langle(2-\lambda)\rho\rangle
                          -m_3/3]\big\}.
\]
The coefficient of $m_3$ in $G-D$ is thus bounded below by
$15/16-5/7+5T/21-65T^2/64$. This concave quadratic exceeds
$3/16$ at $0$ and $1/3$.
Thus
\begin{equation}\label{eq:gdcore}
 G-D\ge \frac3{16}m_3+
 \frac{5T}{7}\left[t(1-t)+\langle(2-\lambda)\rho\rangle-\frac{14}{15}\right].
\end{equation}

\emph{The cutoff moments.} It remains to estimate the two moments of $\rho$. The number
$h=(1-\lambda)t^4+\lambda t^2$ satisfies
$h^2-T(1-\lambda)h-T\lambda=\lambda t^2(h-t^2)\le0$.
Thus $\rho\ge\sqrt h\ge(1-\lambda)t^2+\lambda t$. Integrating the
chord and its cube gives
\[
 \langle(2-\lambda)\rho\rangle+t(1-t)\ge\tfrac53t-\tfrac16t^2,
 \qquad
 \frac{m_3}{T}\ge\frac{t^{-1}+1+t+t^2}{4}\ge\frac89 .
\]
For the last inequality use $t^2\ge t-1/4$ and
$t^{-1}+2t\ge2\sqrt2$. In particular the term $(3/16)m_3$
pays at least $T/6$ of the loss.

Finally \eqref{eq:curvaffine} gives
$h_0(z)/z\le d_*+5T/24-5T^2/48$ on this range.
Together with \eqref{mn:balance}, \eqref{eq:gdcore} and $Z\le13T^2/60$, this bounds
the normalized gap below by
\[
 \eta J-d_*+
 T\left[-\frac{17}{24}+\frac{25}{21}t-\frac5{42}t^2
                              -\frac{11}{20}t^4\right].
\]
Here $t<4/5$ and $t^3<4/9$ since $t^4\le1/3$, so the last three
terms in brackets sum to at least
$(25/21-2/21-11/45)t\ge17t/20$.
The total loss is at most
$(17/24)t^4(1-6t/5)$. Its maximum on $t\ge0$ is
$(17/120)(2/3)^4<1/30$, attained at $t=2/3$.
Since $\eta J-d_*>1/30$, this lower bound for the normalized gap is strictly positive uniformly over $0<T\le1/3$.
\end{proof}

\subsection*{The four-dimensional axis}\label{sec:assequattro}
When $T\ge1/3$, the radii for which the reflected point lies beyond
the axis already supply enough flux. We retain this whole region.
Its length in $x=r^2$ is $\chi=\rho_T(\lambda)^2$, given by
\[
 \chi^2=T[(1-\lambda)\chi+\lambda].
\]
The orbital denominator is convex in $T$, whereas the latitude density
is concave. This permits a comparison at the two endpoint latitudes.

\begin{lemma}\label{lem:outfour}
Let $1/3\le T\le1$, $v=B(1-T)$, $a=B(1-Tq)$, and
$x^2\le q\le x<1$. Put $t=(3T-1)/2$ and
\[
 A(T)=\frac78+\frac7{40}(1-T),\qquad
 R(T)=(1-t)\frac{17}{16}+t\frac74 .
\]
The outgoing density satisfies
\[
 (1-x)w(a)K_{\sqrt{x}}(v,a)\ge w(B)A(T)/R(T).
\]
\end{lemma}
\begin{proof}
Set $h=4x/(1+x)^2$. In \eqref{eq:orbitaltwo} the squared orbital
arguments are $hva$ and $h(1-v)(1-a)$. Since $a\ge v$ and $h\le1$,
the orbital bound gives
\[
 \frac{\mathcal W K_{\sqrt x}(v,a)}{\sqrt{1+x}}
 \ge 1+\frac h{32}\big[20v^2-(1-v)^2-3(1-v)^4\big]
 \ge 1+\frac h{32}(14v-4)\ge A(T).
\]
The second inequality discards $v^2+3v^3(4-v)$. For the last, use
$v\ge2(1-T)/5$ and $28(1-T)/5-4<0$ on this interval.

Write $g(s)=\sqrt{1-s}(1+Bs/C)$, so that $w(a)/w(B)=g(Tq)$.
This function is concave on $[0,1]$. Since $B/C>7/10$,
$g(1/3)\ge(37/30)\sqrt{2/3}>1$ and $g(q/3)\ge1$. At the other
endpoint $g(q)\ge g_1(x)$, where
\[
 g_1(x)=\min\{1,\sqrt{1-x}(1+2x/3)\}.
\]
Indeed $q\in[0,x]$ and $B/C>2/3$. Thus $g(Tq)\ge(1-t)+t g_1(x)$.

The orbital denominator from \eqref{eq:orbitaltwo} can be written
\[
 \frac{\mathcal W}{1-x}
 =(1-x)\big[(1-x)^2+16BCx-8B(1-2B)Tx(1+q)\big]
       +16B^2T^2\frac{x(x-q)(1-xq)}{1-x}.
\]
Both terms decrease with $q$. Their coefficients of degrees $0,1,2$
in $T$ increase with $B$ for $2/5\le B\le5/12$.
Substitution of $q=x^2$ and $B=5/12$ therefore gives the upper bound
\[
 P_x(T)=(1-x)\left[(1-x)^2+\frac{35x}{9}
               -\frac59 Tx(1+x^2)\right]+\frac{25}{9}T^2x^2(1-x^3).
\]
The endpoint bounds, verified in Appendix~\ref{app:scalar}, are
\begin{equation}\label{axis:endpoint}
 P_x(1/3)\le\tfrac{17}{16}\sqrt{1+x},\qquad
 P_x(1)\le\tfrac74\sqrt{1+x}\,g_1(x).
\end{equation}
Convexity of $P_x$ now bounds
$P_x(T)/\sqrt{1+x}$ by $(1-t)(17/16)+t(7/4)g_1(x)$.
Since $g_1\le1$ and $17/16<7/4$, this is at most
$R(T)[(1-t)+tg_1(x)]\le R(T)g(Tq)$, as required.
\end{proof}

\begin{proposition}\label{prop:asse-fine}
For $v=B-BT$ and $1/3\le T\le1$, $\mathcal B_0(BT)>h_0(BT)$.
\end{proposition}
\begin{proof}
The change of variables in \eqref{eq:pairbridge} gives the density
$(1-x)[f_r(v,a)-rf_r(v,a')]$ in $dx\,d\lambda$.
It is nonnegative, and its second term vanishes for $x<\chi$.
Lemma~\ref{lem:outfour} therefore yields
\[
 \mathcal B_0(z)/z\ge w(B) A(T)\langle\chi\rangle/R(T).
\]
For fixed $\lambda$, the positive root
$\chi=pT+\sqrt{p^2T^2+\lambda T}$, $p=(1-\lambda)/2$, increases
and is concave in $T$: the second derivative of the radical is
$-\lambda^2/[4(p^2T^2+\lambda T)^{3/2}]$.
Since $A$ is positive, affine and decreasing,
$A(T)\langle\chi\rangle$ is concave. By \eqref{eq:curvaffine}
it suffices to compare this quantity times $w(B)$ with
$R(T)(3/5+T/10)$, a convex quadratic, at the two endpoints.

At $T=1$ all directions are outgoing and $\chi=1$. The lower bound
there is $w(B)/2>45/64>7/10$.
At $T=1/3$ the root can be written
\[
 \chi=\frac{1-\lambda}{6}+
  \left|\left(\frac{\sqrt{2\lambda}}3,\frac{1+\lambda}{6}\right)\right|.
\]
The mean of the norm is at least the norm of the mean. Hence
$\langle\chi\rangle\ge(3+\sqrt{209})/36>29/60$.
Here $A/R=14/15$ and $w(B)>45/32$, so the lower bound is greater
than $(45/32)(14/15)(29/60)=203/320>19/30$, the required value.
\end{proof}

\begin{proof}[Proof of Theorem~\ref{th:main}]
Proposition~\ref{prop:plus} gives the strict critical inequality on the entire side $v>B$ and a positive lower limit at the link. On the other side Lemma~\ref{prop:otherfifth} gives both conclusions for $0<T\le1/3$, and Proposition~\ref{prop:asse-fine} covers $1/3\le T\le1$. This proves the hypotheses of Proposition~\ref{prop:transfer} on the full closed meridian. The sign of the quartic leaves and the comparison \eqref{eq:cal} follow for every $0<s<s_*$.
\end{proof}

\appendix
\medskip\section{Two one-variable estimates}\label{app:scalar}

\subsection*{The radial integral on the three-dimensional side}
\begin{lemma}\label{lem:radialplus}
For the two radial integrals in \eqref{eq:plusradial} and $1/5\le T\le1$,
\[
 (1+2T/5)(I_1+I_2)\ge51/100.
\]
\end{lemma}
\begin{proof}
Integrate first with measure
$W_T(x)(1+x)^{-2}\,dx$, of mass $N$. If $U$ denotes the integral of
$D_T/(1+x)^2-1$ with respect to this measure, Jensen gives
$I_1\ge N/(1+U/N)$. Substitution of $y=1+x$ gives
\[
 2N=1-(1-\log2)T,\qquad
 U/T=\tfrac67\log2-\tfrac12+\tfrac27T-\tfrac3{16}T^2 .
\]
Using $2/3<\log 2<7/10$ and $T^2\ge2T-1$ yields
$N\ge(1-T/3)/2$, $U/T\le3(1-T/3)/10\le3N/5$.

For $I_2$ set $p=1+2T$ and $d=pT$. Then $D_T(x)\le(1+px)^2$:
even with $6/7$ increased to one the difference is
$Tx[3+(4-2T)x-(1+3T)x^2]\ge0$, by concavity in $x$ of the bracket.
Under $x=Ty/[1+d(1-y)]$, the measure $dx/(1+px)^2$ is $T\,dy/(1+d)$.
Apply Jensen to the two resulting denominators with probability densities
$3\sqrt y/2$ and $2y$. Their means of $1-y$ are $2/5$ and $1/3$.
Consequently
\[
 \frac{1+d}{T}I_2\ge\frac{2\sqrt T/3}{\sqrt{1+2d/5}}
              +\frac{T/2}{1+d/3}
 \ge\frac{2\sqrt T/3}{1+d/5}+\frac{T/2}{1+T}.
\]
This last bound is at least $(1+d)/[3(1+T)]$. Indeed
$\sqrt T\ge(3+7T)/10$ on $[1/5,1]$, by
$(1-\sqrt T)(7\sqrt T-3)\ge0$. Inserting this affine minorant and
clearing the denominators leaves $(1-T)(8T^3+10T^2+19T-4)\ge0$, since the cubic is increasing and positive at $1/5$.

Combining the bounds gives
\[
 (1+2T/5)(I_1+I_2)\ge
 (1+2T/5)\left(\frac{1-T/3}{2(1+3T/5)}+\frac{T}{3(1+T)}\right)
 =\frac12+\frac{T(5+T^2)}{15(1+T)(5+3T)} .
\]
The last term increases, since $(1+T)(5+3T)/T=8+5/T+3T$ decreases.
At $T=1/5$ the right-hand side equals $51/100$.
\end{proof}

\subsection*{Endpoint bounds on the four-dimensional side}
\begin{proof}[Proof of \eqref{axis:endpoint}]
For the first endpoint, expand $P_x(1/3)$ and use $0\le x\le1$:
\[
 P_x(1/3)=1+\tfrac{19}{27}x-\tfrac{32}{81}x^2
      -\tfrac{32}{27}x^3+\tfrac5{27}x^4-\tfrac{25}{81}x^5
 \le 1+\tfrac34x-\tfrac38x^2-x^3 .
\]
The tangent to $3x^2/8+x^3$ at $x=1/4$ is $3x/8-7/128$.
Consequently $P_x(1/3)\le1+7/128+3x/8
\le(17/16)(1+2x/5)\le(17/16)\sqrt{1+x}$.

At the axis the denominator is a square:
$P_x(1)=(1-x)p(x)^2$, $p(x)=1+2x/3+5x^2/3$.
For $x\le1/2$ we have $p(x)\le1+3x/2\le(1+x)^{3/2}$ and
$g_1(x)\ge\sqrt{1-x^2}$, since $1+2x/3\ge\sqrt{1+x}$.
Thus the ratio $P_x(1)/(\sqrt{1+x}\,g_1(x))$ is at most
$\sqrt{1-x}(1+x)^2\le9/(4\sqrt2)<7/4$, by monotonicity up to $1/2$.
For $x\ge1/2$ we have $g_1(x)=\sqrt{1-x}(1+2x/3)$, since this product decreases there and is below one at $x=1/2$. Using
$\sqrt{1+x}\ge1+2x/5$ gives
$\sqrt{1+x}(1+2x/3)\ge(1+x/2)^2$.
Moreover $p(x)\le(7/15)(1+4x)(1+x/2)$, the difference being
$(2x-1)(16-11x)/30$. The ratio is now at most
\[
 (7/15)^2\sqrt{1-x}(1+4x)^2\le8(7/15)^2<7/4,
\]
where the maximum of the radical times the square occurs at $x=3/4$.
This proves \eqref{axis:endpoint}.
\end{proof}

\bigskip
\section*{AI methodology}
The author began this project several years ago. His earlier manuscript gave
a radial criterion for minimality and local estimates at the Lawson link.
What remained was to prove the sign of the curvature on each side of the
cone. Starting from that manuscript, the author used OpenAI's ChatGPT to
explore radial inversion and cumulative inequalities. ChatGPT helped apply
Cabr\'e's calibration to the quartic leaves, derive the reflection inequalities
at reciprocal radii in Proposition~\ref{prop:quarticmirror}, and establish
the continuity at the link needed to pass to positive orders in
Proposition~\ref{prop:transfer}.

For the quadratic curvature, the author proposed a chord bound for the
normalized excess over $d_*$ and identified the logarithmic estimate needed
at $T=1$. In Lemma~\ref{lem:four-curvature}, ChatGPT applied this idea to
$I(T)$ and used the bounds on the median to estimate $I(1)$.

In the initial work on the angular kernel, the author suggested using
H\"older's Inequality to bound the moment of order $7/2$ by integer moments.
ChatGPT supplied the substitution in the semicircle integral that yields the
quadratic term in~\eqref{eq:orbitaltwo}. The author also proposed deriving
the diagonal kernel estimate from the Poisson semigroup and Cauchy--Schwarz.
ChatGPT helped extend it to distinct latitudes in~\eqref{eq:poissonorbit}
and evaluate the mass of the radial measure at the link.

ChatGPT's contribution to the proof of Theorem~\ref{th:main} was decisive.
It derived the analytic estimates and carried out the arithmetic comparisons
that establish the flux inequalities on both sides of the cone. These include
the radial bound in Lemma~\ref{lem:radialplus}, the comparison of gain and
loss at a common cutoff in Lemma~\ref{prop:otherfifth}, and the endpoint
estimates and lower bound for $\langle\chi\rangle$ at $T=1/3$ used in
Proposition~\ref{prop:asse-fine}.

The author takes responsibility for the mathematical content and the references.

\end{document}